\documentclass[12pt,reqno]{amsart}
\usepackage{amsmath,amssymb}
\usepackage[dvips]{graphicx,color}
\usepackage{latexsym,amssymb}
\usepackage{mathrsfs}
\usepackage[abbrev]{amsrefs}
\usepackage{comment}
\usepackage[dvipdfmx]{xcolor}
\usepackage{latexsym}
\def\qed{\hfill $\Box$}

\newtheorem{thm}{Theorem}[section]

\newtheorem{prop}[thm]{Proposition}
\newtheorem{lem}[thm]{Lemma}

\newcommand{\ep}{\varepsilon}

\def\<{\langle }

\newcommand{\supp}{\text{ supp }}

\numberwithin{equation}{section}
\numberwithin{thm}{section}

\begin{document}

\begin{center}\large \bf 
On the ill-posedness
for the full system of compressible Navier-Stokes equations
\end{center}

\footnote[0]
{


E-mail: 
\, $^*$motofumi.aoki.q2@dc.tohoku.ac.jp
$^{**}$t-iwabuchi@tohoku.ac.jp
}
\vskip5mm

\begin{center}
{\sf 
Motofumi Aoki$^{*}$
Tsukasa Iwabuchi$^{**}$
}

\vskip2mm

Mathematical Institute, 
Tohoku University\\
Sendai 980-8578 Japan

\end{center}

\vskip5mm

\begin{center}
\begin{minipage}{135mm}
\footnotesize
{\sc Abstract.}
We consider the Cauchy problem for compressible Navier--Stokes equations of the ideal gas in the three-dimensional spaces. 
It is known that the Cauchy problem in the scaling critical spaces of the homogeneous Besov spaces $\dot B^{\frac{3}{p}}_{p,1}\times\dot B^{-1+\frac{3}{p}}_{p,1}\times\dot B^{-2+\frac{3}{p}}_{p,1}$ is uniquely solvable for all $1 < p<3$ and is ill-posed for all $p>3$. 
However, it is an open problem whether or not it is well-posed 
in the case when $p=3$. 
In this paper, we prove that for the case $p=3$ the Cauchy problem is ill-posed by constructing a sequence of initial data, which 
shows that the solution map is discontinuous. 
\end{minipage}
\end{center}

\section{Introduction}
We consider the Cauchy problem for 
the full-system of compressible Navier-Stokes equations in $\mathbb R^3 $:
\begin{equation}\tag{cNS} \label{cNS4}
\left\{
	\begin{split}
	&\partial_t \tilde{\rho} + {\rm div \,} (\tilde{\rho} u) =0, 
	&t > 0, \, x \in \mathbb R^3 , \\
	&\partial_t (\tilde{\rho} u)  + {\rm div}\,( \tilde{\rho} u \otimes u)  +  \nabla  (\tilde{\rho} \theta) =  \mu \Delta u + (\mu +\lambda) \nabla {\rm div \,} u ,
	&t > 0, \, x \in \mathbb R^3 , \\
	&\partial_t (\tilde{\rho} \theta) + {\rm div}\, (\tilde{\rho} \theta u) +\tilde{\rho} \theta \, {\rm div \,} u - \kappa \Delta \theta = 2\mu |D(u)|^2 + \lambda |{\rm div \,}u|^2,
	&t > 0, \, x \in \mathbb R^3 , \\
&( \tilde{\rho}, u, \theta)\big|_{t=0} = (\tilde{\rho}_0, u_0, \theta_0),
&x \in \mathbb R^3. 
	\end{split}
\right.	
	\end{equation}
    where $\mu $ and $\lambda$ are the Lam\'e constants satisfying $\mu >0, \, 2\mu + \lambda >0$ and $ \kappa$ is a positive constant.
    The unknown functions are 
    $\tilde{\rho} = \tilde{\rho}(t,x)$, $u = (u_1(t,x), \cdots, u_d (t,x))$ and $\theta =\theta(t,x)$ which are the density, velocity and temperature of the fluid, and $\tilde{\rho}\theta$ denotes the pressure given by the equation of the ideal gas. 
We note that the other physical constants except for $\mu, \lambda , \kappa$ are taken by $1$ for simplicity. 
    The deformation tensor is given by 
    \[
        D(u) = \frac{\nabla u +(\nabla u)^{T}}{2}, 
    \]
    where $\nabla u =(\partial_{x_i}u_j)_{i,j = 1}^3$, $|A|^2$ denotes the trace of $AA^T$ for a matreix $A$ and its transpose $A^T$.
    We assume that the density of the fluid satisfies the non-vacuum condition. 
    \[
        \inf_{t\geq 0, x \in \mathbb R^3}\tilde \rho (t,x) > 0. 
    \]
More precisely, we write the density 
$\tilde \rho = 1 + \rho$ and rewrite \eqref{cNS4} as follows. 
    \begin{equation}
        \label{cNS5} 
        \left\{
            \begin{split}
        &\partial_t \rho + {\rm div \,} u + 
        {\rm div \,}(\rho u) =0, 
        &t > 0, \, x \in \mathbb R^3 , \\
        &\partial_t u  - \mathcal{L} u + (1+\rho)
        (u \cdot \nabla) u + \nabla ((1+ \rho)\theta) = -\rho\partial_t u,
        &t > 0, \, x \in \mathbb R^3 , \\
        &\partial_t \theta +(1+\rho)(u \cdot \nabla) \theta + (1+\rho) \theta \, {\rm div \,} u - \kappa \Delta \theta
        \\
        &\hspace{44pt}
        = 2 \mu |D(u)|^2 + \lambda |{\rm div \,}u|^2 - \rho \partial_t \theta, 
        &t > 0, \, x \in \mathbb R^3 , \\
        &( \rho, u, \theta)\big|_{t=0} = (\rho_0, u_0, \theta_0),
        &x \in \mathbb R^3,  
    \end{split}
    \right.	
    \end{equation}
    where we set $\mathcal{L} = \mu \Delta +(\lambda + \mu) \nabla {\rm div \,}$. 
    Then we consider the existence of the strong solution under the uniformly parabolic condition for the velocity and the temperature equations.
    We define the Besov spaces as follows:
    \vskip3mm

    \noindent
    {\bf Definition}
    {\rm \,(The Besov spaces)}.
    Let $\{\phi_j\}_{j\in \mathbb{Z}} \subset \mathcal S(\mathbb R^3)$ be such that 
    \[
    \begin{aligned}
      &\supp  \widehat{\phi_0} \subset
      \{\xi \in \mathbb{R}^3 |\,2^{-1} \le |\xi| \le 2 \},
      \hspace{33pt}
      \phi_j (x) = 2^{3j} \phi_0 (2^j x),
      \\
      &
      \sum_{j\in \mathbb{Z}} \widehat{\phi_j}(\xi) =1
      \hspace{11pt} {\rm for} \hspace{11pt}
      \xi \in \mathbb R^3 \backslash \{0\}.
    \end{aligned}
    \]
   For $s\in \mathbb R$, $1 \le p,q \le \infty$, we define that
   \[
    \dot B^s_{p,q} (\mathbb R^3) :=
    \bigg\{
      f \in \mathcal S'/\mathcal P \bigg|
      \|f\|_{\dot B^s_{p,q}} := 
      \big(\sum_{j\in \mathbb Z} 2^{sjq} 
      \|\phi_j * f\|_p^q \big)^{\frac{1}{q}} < \infty
    \bigg\},
   \]
where $\mathcal P$ denotes the set of all polynomials. 

\vskip3mm 

   Let us recall the existing results about the well-posedness and the ill-posedness of the Cauchy problem for the compressible Navier--Stokes equations. 
   Furthermore, we will compare with the known results of energy conseravtion law.
   \color{black}
\\
   
   \vskip3mm
   \noindent
   {\bf A. Strong solution.}

The scale invariance for the equations \eqref{cNS4} can be understood as follows. 
   For any $\alpha >0$, the scaling  transformation of $(\rho, u, \theta)$ given by 
   \[
   \begin{cases}
       \rho_{\alpha}(t,x) = \rho(\alpha^2 t, \alpha x), \\
       u_{\alpha}(t,x) = \alpha u(\alpha^2 t, \alpha x), \\
       \theta_{\alpha}(t,x) = \alpha^2 \theta(\alpha^2 t, \alpha x)
   \end{cases}
   \]
   maintains the equations \eqref{cNS4} and the norms of the initial data are invariant with respect to $\alpha>0$ in function spaces such as $\dot H^{\frac{d}{p}}_p \times \dot H^{-1+\frac{d}{p}}_p \times \dot H^{-2+\frac{d}{p}}_p$ which are called scaling invariant spaces.

   Danchin \cite{Da-011}, \cite{Da-012} first considered the barotropic equations conrresponding to \eqref{cNS5} $(\theta \equiv 0)$ in the scaling critical Besov spaces such that
   \[
   \begin{aligned}
       \rho \in C([0,T); \dot B^{\frac{d}{p}}_{p,1}(\mathbb{R}^d)), \, 
       u \in C([0,T); \dot B^{{-1+\frac{d}{p}}}_{p,1}(\mathbb{R}^d)),
   \end{aligned}
   \]
   with $1 < p <d$. 
   Later Haspot \cite{Ha-11} generalized the local well-posedness result to $(\rho_0, u_0) \in \dot B^{\frac{d}{p}}_{p,1} (\mathbb{R}^d)\times \dot B^{-1+\frac{d}{p}}_{p,1} (\mathbb{R}^d)\, (1 < p <2d)$.
   Furthermore, Chikami--Danchin \cite{ChDa-15} discuss the unique solvability with $(\rho_0 -1, u_0, E_0) \in (\dot B^{\frac{d}{p}}_{p,1}, \dot B^{-1+\frac{d}{p}}_{p,1}, \dot B^{-2+\frac{d}{p}}_{p,1})$ $(1 < p < 2d)$, where $E_0:= \displaystyle \frac{\rho_0 u^2_0}{2} + e(\rho_0, \theta_0)$ is the total energy, and that with $(\rho_0 -1, u_0, \theta_0) \in (\dot B^{\frac{d}{p}}_{p,1}, \dot B^{-1+\frac{d}{p}}_{p,1}, \dot B^{-2+\frac{d}{p}}_{p,1})$ $(1 < p < d)$.

   On the other hand, there are several studies about the ill-posedness of the motion of the ideal gas. 
Chen--Miao--Zhang \cite{ChMiZh-15} proved the ill-posedness in the case when $d=3$ for the initial data $(\rho_0, u_0, \theta_0) \in (\dot B^{\frac{3}{p}}_{p,1}, \dot B^{-1+\frac{3}{p}}_{p,1}, \dot B^{-2+\frac{3}{p}}_{p,1}) $ with $p>3$.
They also show the ill-posedness of the Cauchy problem for the barotropic case with $p >2d$. 
Recently, Iwabuchi and Ogawa \cite{IwOg-21} proved the ill-posedness in the case when $d=2$ for the initial data 
\begin{equation}
    \label{110501}
    (\rho_0, u_0, \theta_0) \in (\dot B^{\frac{2}{p}}_{p,q}, \dot B^{-1+\frac{2}{p}}_{p,q}, \dot B^{-2+\frac{2}{p}}_{p,q})   
\end{equation}
 with $1 \le p \le \infty$ and $1 \le q < \infty$.
They~\cite{IwOg-22} also show the ill-posedness for the barotropic equations for $p=2d$.

    \vskip3mm
   \noindent
   {\bf B. Weak solutions and the energy conservation law.}

    Lions \cite{Li-98} first showed the existence of weak solutions for the generalized isentropic Navier--Stokes equations on the bounded domain. 
    Later Feireisl \cite{Fe-04} showed the existence of weak solutions for compressible Navier--Stokes equations in $\Omega$, where $\Omega$ is a smooth bounded domain.

There are many studies about the energy conservation law for the compressible Navier--Stokes equations.
Yu \cite{Yu-17} proved the energy identity for isentropic Navier-Stokes equations 
if the velocity belongs to $L^p(0,T; L^q(\mathbb T^3))$ with $ 2/p + 2/q \le 5/6$ and $q \ge 6$. 
Akramov--D\polhk ebiec--Skipper--Wiedemann \cite{AkDeSkWi-20} showed the energy identity when   
$u \in B^{\alpha}_{3, \infty}((0,T) \times \mathbb T^3) \cap L^2(0,T;W^{1,2})$, $\rho, \rho u \in B^{\beta}_{3, \infty}((0,T) \times \mathbb T^3)$ with $\alpha + 2\beta >1,\,  2 \alpha + \beta >1,\, 0 \le \alpha, \beta \le 1$ and 
the pressure is continuous with respect to the density. 
Furthermore, Nguyen--Nguyen--Tang~\cite{NgNgTa-19} 
considered weak solutions such that 
$u \in L^p(0,T;L^q), \, \rho \in L^{\infty}(0,T; B^{\frac{1}{2}}_{\frac{12}{5},\infty}(\mathbb T^3))$ with $2/p + 2/q=1, \, q \ge 4$, 
and the pressure is a $C^2$ function with respect to the density,  
to prove energy identity. 
Recently the authors \cite{AoIw-22} showed the energy conservation law for the full system of the compressible Navier--Stokes equations and the barotropic equations for two and three dimensional case if the density satisfies $\rho \in L^{\infty}(0,T;L^{\infty})$ and the velocity satisfies $u \in L^{\infty}(0,T;L^2)$ 
 and in the case when $d=3$, 
\color{black}
    \begin{equation}
    \label{091505}
    \begin{split}
    \displaystyle 
    & 
    u \in L^p(0,T; L^q), 
    \quad \frac{1}{p} + \frac{3}{q} = 1,
    \quad 3 \le q \le 4,
    \\
    \displaystyle
    \text{or } \quad 
    & u \in L^p(0,T; L^q),
    \quad \frac{2}{p} + \frac{2}{q}= 1, 
    \quad q \ge 4 
    \end{split}
    \end{equation}
    Moreover, assuming the existence of a weak solution for the ideal gas equations \eqref{cNS4}, we also obtain that the weak solution satisfies the energy conservation law under the above condition.

    \vskip3mm
    In this paper, we study whether the Cauchy problem for the equations \eqref{cNS5} with the initial data $(\rho_0, u_0, \theta_0) \in (\dot B^{1}_{3,1}, \dot B^{0}_{3,1}, \dot B^{-1}_{3,1})$ is ill-posed.
    Furthermore, we discuss the difference between the class of the Cauchy problem and the class of energy conservation law. 

    \color{black}
    Looking at the previous studies on the barotropic case, the additional condition of \eqref{091505} is a reasonable result comparing with the result of the Cauchy problem.
    In fact, the barotropic equations with data  $(\rho_0, u_0) \in  (\dot B^{1}_{3,1}, \dot B^{0}_{3,1})$ are uniquely solvable (see Haspot \cite{Ha-11}).
    If $\rho \in L^{\infty}(0,T;L^{\infty})$ and $u \in L^{\infty}(0,T;L^3)$, the bartropic equations satisfies energy conservation law (see \cite {AoIw-22}).
    \color{black}
    Comparing with above conditions,
    the class of the Cauchy problem includes the class of energy conservation law as follows;
    \begin{equation}
      \label{111002}
      \rho \in C([0,T); \dot B^{1}_{3,1}) \subset L^{\infty}(0,T; L^{\infty}), 
      \quad 
      u \in C([0,T); \dot B^{0}_{3,1}) 
      \subset L^{\infty}(0,T; L^3).
    \end{equation} 
    In the case of a two-dimensional ideal gas, however, the relation between the class of the Cauchy problem and the class of energy conservation law seems conflicting.
    \color{black}
    It is known that the following inclusion relations hold:
    \[
      C([0,T) ;\dot B^{1}_{2,1}) \subset L^{\infty}(0,T;L^{\infty}), \,
      C([0,T) ;\dot B^{0}_{2,1}) \subset L^{\infty}(0,T;L^2).
    \]
    $\rho \in L^{\infty}(0,T;L^{\infty}),\,u \in L^{\infty}(0,T;L^2)$ are sufficient conditions that the ideal gas satisfies the energy conservation law (\cite{AoIw-22}).
    On the other hand, 
    $(\rho_0, u_0, \theta_0) \in (\dot B^{1}_{2,1}, \dot B^{0}_{2,1}, \dot B^{-1}_{2,1})$ is one of the initial data for which Iwabuchi--Ogawa \cite{IwOg-21} showed that the Cauchy problem is ill-posed.
    There is indeed no relation between $L^1$ and $\dot B^{-1}_{2,1}$. 
    However, except for this point, these results imply that the equations \eqref{cNS4} has different aspests under similar conditions: a positive result that the equations \eqref{cNS4} satisfy the energy conservation law and a negative result that the Cauchy problem of \eqref{cNS5} is ill-posed. 

    \color{black}


\color{black}
\begin{thm}
The Cauchy problem \eqref{cNS5} is ill-posed in $\dot B^1_{3,1} \times \dot B^0_{3,1} \times \dot B^{-1}_{3,1}$. 
More precisely, there exist a sequence of initial data
$\{u_{0,N} \}_N\subset \mathcal S(\mathbb R^3)$, positive time $\{T_N\}_N$ with $T_N \to 0$ as $N \to \infty$, and a sequence of corresponding smooth solutions $\{(\rho_N, u_N, \theta_N)\}_N$ in the time interval $[0,T_N]$ with the data $(0,u_{0,N}, 0)$ such that
\[
\lim_{N \to \infty} \|u_{0,N}\|_{\dot B^0_{3,1}} = 0, \quad
\lim_{N \to \infty} \|\theta_n (T_N)\|_{\dot B^{-1}_{3,1}} = \infty.
\]
\end{thm}

\noindent
{\bf Remark.}
\begin{enumerate}
  \item 
Chikami--Danchin \cite{ChDa-15} discussed the unique solovability in the case when $1 \le p<d$, and 
  Chen--Miao--Zhang \cite{ChMiZh-15} proved the ill-posedness result in the case when $p>d$. 
  The threshold case is settled by showing the ill-posedness.
  \item 
If we restrict to the case when $p=d$,
Theorem 1.1 also implies a similar relation between the class of the Cauchy problem and the class of energy conservation law like the two-dimensional ideal gas case.
  In fact, we prove the ill-posedness for the equation \eqref{cNS5} with the initial data $(\dot B^1_{3,1}, \dot B^{0}_{3,1}, \dot B^{-1}_{3,1})$.
  On the other hand, the equations \eqref{cNS4} in $\rho \in L^{\infty}(0,T;L^{\infty})$ and $u \in L^{\infty}(0,T;L^3)$ satisfy energy conservation law. 
  Then we also have the corresponding relations for the density and the velocity like \eqref{111002}. 
  
\end{enumerate}
\color{black}

\section{preliminary}
We introduce the modulation spaces, the bilinear estimate and the maximal regularity. 

\vskip3mm
\noindent
{\bf Definition} 
  {\rm \,(The modulation spaces~\cite{Fe-83})}. 
We say that $f \in \mathcal{S}'(\mathbb{R}^3)$ 
belongs to  $M_{3,1}(\mathbb{R}^3)$ if  
  \[
    \|f\|_{M_{3,1}} := 
    \sum_{k \in \mathbb Z^3}
    \|\mathcal F^{-1}[\chi(\xi-k)\hat{f}]\|_3
    < \infty,
  \]
  where $\chi \in C^{\infty}_0(\mathbb{R}^3)$ satisfies 
  \[
  \begin{aligned}
   & \supp \chi \subset \{ \xi \in \mathbb R^3| \xi_j \in [0, 0 +2), \,j =1,2,3\},\\
   & \sum_{k \in \mathbb{Z}^3} \chi (\xi -k) =1 
   \quad \text{for all } \xi \in \mathbb R^3. 
  \end{aligned}
  \]

\begin{lem}
There exists $C>0$ such that the following bilinear estimate holds.
\begin{equation}
    \label{102403}
    \|f g\|_{M_{3,1}}
    \le
    C\|f\|_{M_{3,1}}
    \|g\|_{M_{3,1}}.
\end{equation}
\end{lem}

\noindent
The above estimate \eqref{102403} is obtained by the following bilinear estimate and the embedding theorem (see e.g. \cite{To-04}).
\[
  \|fg\|_{M_{3,1}} 
  \le
  C\|f\|_{M_{\infty,1}}  \|g\|_{M_{3,1}}, 
  \hspace{22pt}
  \|f\|_{M_{\infty,1}}  
  \le 
  C  \|f\|_{M_{3,1}}.
\]

\noindent
{\bf Proposition 2.2 (Maximal regularity~\cite{IwOg-21})}.
Let $f \in M_{3,1}, u \in L^1(0,T;M_{3,1})$.
There exists $C>0$ independent of $f,u$ such that for any $T >0$,
\begin{equation}
  \label{102404}
      \|\mathcal{L} e^{t\mathcal{L}}f\|_{L^1(0,T;M_{3,1})}
    \le
    C \sqrt{\log (e+T)} \|f\|_{M_{3,1}},
\end{equation}
\begin{equation}
  \label{102405}
    \big\|
    \int^t_0
     \mathcal{L} e^{(t-s)\mathcal{L}}u(s) ds\big\|_{L^1(0,T;M_{3,1})}
     \le 
     C  \sqrt{\log (e+T)} \|u\|_{L^1(0,T;M_{3,1})}.
\end{equation}

\section{proof of Theorem 1.1}

We take an initial data 
such that for $N \in \mathbb{N}$ and $x \in \mathbb{R}^3$, 
\begin{equation}\label{103101}
\begin{aligned}
  &\rho_0(x) =0, 
  \hspace{11pt}
    u_{0,N}(x) = (u_{0,N,1}(x),0,0),
    \hspace{11pt} 
    \theta_0(x) =0,
    \vspace{5pt}
    \\
   & \begin{cases}
      \displaystyle
      u_{0,N,1}(x) \hspace{-11pt}
      &= 
      \displaystyle
      R N^{-\frac{1}{3}} \sum_{- \delta N \leq j \leq 0} 
2^{-2j} \phi_j (x - 2^{|j| +2N}e_1) {\rm sin} (2^Nx_1),
\vspace{5pt}
\\
\displaystyle
   \hspace{33pt}R \hspace{-11pt}
   &= ({\rm log \,} N)^{-1}, 
    \end{cases}
    N = 1, 2, \cdots,
\end{aligned}
\end{equation}
where $e_1 = (1,0,0) $ and $ 0 < \delta < 1$.
We introduce the following integral equations of \eqref{cNS5}. 
\[
\begin{cases}
    \displaystyle
    \rho (t) &
    \displaystyle=- \int^t_0
   {\rm div \,} (u(s)+ \rho(s)u(s))ds, \\
   \displaystyle
   u(t) &
   \displaystyle= e^{t \mathcal{L}} u_{0,N} + 
   \int^t_0 e^{(t-s)\mathcal{L}}
   \big\{
  -(1+\rho(s))(u \cdot \nabla) u(s) - \nabla p(s) - \rho \partial_t u(s) 
   \big\} ds, \\
   \displaystyle
   \theta (t) &
   \displaystyle= 
   \int^t_0 e^{(t-s)\kappa \Delta} 
   \big\{
    -(1 +\rho(s))(u \cdot \nabla)\theta(s) - p(s)\,  {\rm div \,}u(s)
    -\rho(s)\partial_t \theta(s) \\
    &\hspace{55pt}
    +2\mu|D(u(s))|^2 + \lambda ({\rm div \,} u(s))^2 
  \big\}
   ds , 
\end{cases}
\]
\noindent
where $p= (1 +\rho)\theta$.
Let us define the expansion $\{P_k, U_k, \Theta_k\}_{k=0}^{\infty}$ for the initial data $(0,u_{0,N},0)$, 
\[
\begin{cases}
  \displaystyle
    P_0 (t):= 0,
    \quad
    P_1 (t):= - \int^t_0 {\rm div \,} U_1 (s)ds, 
    \\
    \displaystyle
    P_k (t):= 
    -\int^t_0
    \bigg({\rm div \,} U_k (s) + 
    \sum_{k_1+k_2 =k} {\rm div \,}\big(P_{k_1}(s)U_{k_2}(s)\big)\bigg) ds, \qquad k \ge 3,
\end{cases}
\]
\[
\begin{cases}
  \displaystyle
    U_0(t) := 0, 
    \quad
    U_1(t) := e^{t \mathcal{L}} u_{0,N}, \\
    \displaystyle
    U_2(t) := 
    \int^t_0 e^{(t-s)\mathcal{L}}
    \bigg(
      - (U_1 (s)\cdot \nabla) 
    U_1(s)
    -
    \nabla \Theta_2(s) - P_1(s) \partial_t U_1(s) 
    \bigg) ds, 
    \\
    \displaystyle
    U_k(t) := 
    \int^t_0 e^{(t-s)\mathcal{L}}
    \bigg\{
      - \sum_{k_1 +k_2=k} (U_{k_1} (s)\cdot \nabla)U_{k_2}(s)
      \displaystyle
      - \sum_{k_1+k_2+k_3 =k}
      P_{k_3}(s) (U_{k_1}(s) \cdot \nabla) U_{k_2}(s)\\
      \hspace{66pt}
      \displaystyle
      -  \nabla\bigg( \Theta_k (s) +
      \sum_{k_1+k_2 =k} P_{k_1}(s) \Theta_{k_2} (s)
      \bigg)
      - \sum_{k_1+k_2=k} P_{k_1} (s)\partial_t U_{k_2} (s) 
    \bigg\}ds, \quad k \ge 3,
\end{cases}
\]
\[
\begin{cases}
    \Theta_0(t) := 0, \quad \Theta_1(t) :=0,\\
    \displaystyle
    \Theta_2(t) :=
    \int^t_0 
    e^{(t-s)\kappa \Delta}
    \bigg(
      \tilde{D}(U_1)(s) : \tilde{D} 
    (U_1)(s) \bigg) ds,
    \\
    \displaystyle
    \Theta_k(t) :=
    \int^t_0 
    e^{(t-s)\kappa \Delta}
    \bigg\{
    -\sum_{k_1+k_2=k} (U_{k_1}(s)\cdot \nabla) \Theta_{k_2}(s)
    - \sum_{k_1+k_2+k_3 =k} 
    P_{k_3} (s) (U_{k_1}(s)\cdot \nabla) \Theta_{k_2}(s)
    \\
    \displaystyle
    \hspace{99pt}
    -\bigg(
    \sum_{k_1+k_2 =k} \Theta_{k_1}(s) \, {\rm div \,} U_{k_2} (s) +
    \sum_{k_1+k_2+k_3 =k}
    P_{k_3}(s)\Theta_{k_1} (s) 
    {\rm div \,} U_{k_2}(s) \bigg)
    \\
    \displaystyle
    \hspace{99pt}
    -\sum_{k_1+k_2=k} P_{k_1}(s) \partial_t \Theta_{k_2}(s) 
    +\sum_{k_1+k_2=k} \tilde{D}(U_{k_1})(s) : \tilde{D}(U_{k_2})(s)
    \bigg\} ds, \quad k\ge 3, 
\end{cases}
\]

\noindent
where we have used the notation $\tilde{D}(u)(s) : \tilde{D}(v)(s) := 2\mu \, {\rm tr \,} (Du \cdot Dv^{T}) + \lambda ({\rm div \,} u \cdot {\rm div \,}v)$ and all the summations means
\[
  \begin{aligned}
    \sum_{k_1+k_2=k} &a_{k_1} b_{k_2}=
    \sum_{k_1=1}^{k-1} a_{k_1} b_{k-k_1}, \qquad 
    \sum_{k_1+k_2+k_3=k} &a_{k_1} b_{k_2}c_{k_3}=
    \sum_{k_1=1}^{k-1} \sum_{k_2=1}^{k-k_1}a_{k_1} b_{k_2} c_{k-(k_1+k_2)}. 
  \end{aligned}
\]

We first show the upper bound of the initial data and the lower bound of the right-hand side of the temperature term $\Theta_2(t)$.

\begin{prop}
Let $u_{0,N}$ be defined by \eqref{103101}. 
Then, there exists $C>1$, $0 < c< 1$ and sufficient small $\ep_0 >0$ such that for any $N \in \mathbb N$, 
\begin{align}
    \label{110401}
\displaystyle
&\| u_{0,N} \|_{\dot B^0_{3,1}} \le CR \delta^{\frac{1}{3}}, 
\\
    \label{110402}
    \displaystyle
&\bigg\| \int^t_0 e^{(t - s)\mathcal L} | \nabla 
e^{ s\mathcal L} u_{0,N} |^2 ds 
 \bigg\|_{\dot B^{-1}_{3,1}} \ge c\delta R^2 N^{\frac{1}{3}} -
 C(R N^{-\frac{1}{3}})^2 2^{-cN},
\end{align}
for $ t = \varepsilon_0 2^{-2N}$.
\end{prop}

\noindent
{\bf Proof of Proposition 3.1. }
We prove the lower estimate \eqref{110402} 
in more general dimensions 
in the appendix \eqref{081401}.

We show \eqref{110401}. 
Since
\[
  \vspace{-10pt}
\supp \mathcal F[u_{0,N}] \subset 
\{ \xi \in \mathbb R^3
\big| 2^{-(\delta N +1)} \le |\xi - 2^{2N} e_1| \le  2
\}.
\]
We immediately obtain 
\[
  \| u_{0,N} \|_{\dot B^0_{3,1}} 
  \le 
  C\|(\phi_{N-1} + \phi_{N}+\phi_{N+1}) * u_{0,N} \|_3 
  \le
 C\|u_{0,N}\|_3.
 \]
 We set a subset of $(j,k,l)$ such that 
 \[
 \displaystyle J= \{(j,k,l) | - \delta N \le j,k,l \le 0, \, j=k=l\},  
 \quad \displaystyle J^c = 
\{(j,k,l) | - \delta N \le j,k,l \le 0, \,  
  (j,k,l) \not \in J\} ,
\]
and introduce 
\[
\Phi^{\pm}_{N,j}(x) := 
\mathcal F^{-1} [ \widehat \phi_j(\xi \mp 2^{2N} e_1)] (x - 2^{2N+|j|}e_1) .
\]
We then write 
\[
u_{0,N} 
= R N^{-\frac{1}{3}} \sum_{-\delta N \leq j \leq 0} 
2^{-2j} \dfrac{1}{2i} (\Phi_{N,j}^+ - \Phi_{N,j}^-), 
\]
 \[
 \begin{aligned}
  &|u_{0,N}(x)|^3
\\
\le & \dfrac{R^3 N^{-1}}{8}
  \sum^{0}_{-\delta N \leq j \leq 0}
  \sum_{-\delta N \leq k \leq 0}
  \sum_{-\delta N \leq l \leq 0} 2^{-2(j+k +l)}
\Big(|\Phi_{N,j}^+ | + | \Phi_{N,j}^- |\Big)
\Big(|\Phi_{N,k}^+ | + | \Phi_{N,k}^- |\Big)
\Big(|\Phi_{N,l}^+ | + | \Phi_{N,l}^- |\Big)
  \\
  =
  & 
  \sum_{(j,k,l)\in J} + \sum_{(j,k,l)\in J^c} =:A_1+A_2. 
 \end{aligned}
 \]
To prove \eqref{110401}, we first consider the estimate for $A_1$.
Since
\[
\begin{aligned}
  &\int_{\mathbb{R}^3}
\Big(|\Phi_{N,j}^+ | + | \Phi_{N,j}^- |\Big)^3
 dx 
  \le C2^{6j}, 
\end{aligned}
\]
we obtain
\[
\begin{aligned}
    &\|A_1\|_3 \le
    RN^{-\frac{1}{3}} \bigg( \sum_{(j,k,l)\in J}
2^{-2(j+j+j)} \cdot C 2^{6j} 
\bigg)^{\frac{1}{3}}
    \le
    C R  \delta^{\frac{1}{3}}.
\end{aligned}
\]
We next consider the estimate for $A_2$.
In the case when two of $j,k,l$ are same and the other is different. It is sufficient to consider the case when 
$j=l \neq k$. 
We notice that 
the supports of $\Phi_{N,j}^\pm, \Phi_{N,k}^\pm$ 
are far from each other. 
Since for any $m \in \mathbb N$, 
there exists $C > 0$ such that 
\[
\begin{aligned}
  &\int_{\mathbb{R}^3}
\Big(|\Phi_{N,j}^+ | + | \Phi_{N,j}^- |\Big)^2
\Big(|\Phi_{N,k}^+ | + | \Phi_{N,k}^- |\Big)
 dx 
  \le C 2^{6j+3k} 2^{-mN} , 
\end{aligned}
\]
we have 
  \begin{equation}
  \begin{aligned}
    \label{110702}  
&   \dfrac{R^3 N^{-1}}{8} \sum_{-\delta N \le j \le 0}
    \sum_{j\neq k}
    2^{-2(j+k+j)}
\int_{\mathbb R^3} 
\Big(|\Phi_{N,j}^+ | + | \Phi_{N,j}^- |\Big)^2
\Big(|\Phi_{N,k}^+ | + | \Phi_{N,k}^- |\Big)
    dx 
\\
\le&
 C R^3 N^{-1}  \sum_{-\delta N \le j \le 0}
    \sum_{j\neq k}
    2^{-2(2j+k)+6j +3k}2^{-mN}
\leq C R^3 N^{-1}2^{ - mN}. 
  \end{aligned}
  \end{equation}

If the three of $j,k,l$ are different each other, a similar estimate implies the same bound, 
which yields \eqref{110401}.
\hfill $\Box$

\vskip3mm 

\color{black}
We finally prove the upper bound of $(P_k, U_k, \Theta_k)$.

\begin{prop}
    \label{103103}
    Let $u_{0,N}$ be defined by \eqref{103101}.
    Then there exist $C_0 >0$ such that for all $t \le \varepsilon_0 2^{-2N}$, 
    \begin{equation}
        \label{103102}
        (t 2^N)^{-1} \|P_k(t)\|_{M_{3,1}} + 
        \|U_k (t)\|_{M_{3,1}} + 
        2^{-N} \|\Theta_k(t)\|_{M_{3,1}} 
        \le C_0^k t^{k-1} 2^{(k-1)N} R^k.
    \end{equation}
    \end{prop}

\noindent
{\bf Proof of Proposition~\ref{103103}. }
Since $\supp \mathcal{F} [u_{0,N}]$ have the frequency around $2^N$, 
\begin{equation}
\begin{aligned}
    \label{102801}
    &\supp \mathcal{F}[P_k], 
    \supp \mathcal{F}[U_k], 
    \supp \mathcal{F}[\Theta_k]
    \\
    & \hspace{22pt}
    \subset 
    \{ \xi \in \mathbb{R}^3|
    |\xi| \le k 2^{N+1}
    \},
    \qquad k=1,2,3, \cdots.
\end{aligned}
\end{equation}
We calculate directly in the case when $k=1$.
Here $C>0$ exists such that 
\[
\begin{cases}
    \displaystyle
    \|U_1(t)\|_{M_{3,1}} = 
    \|e^{t \mathcal{L}}u_{0,N}\|_{M_{3,1}}
    \le \|u_{0,N}\|_3
    \le  CR, 
    \vspace*{3pt}
    \\
    \displaystyle
    \|\nabla U_1(t)\|_{M_{3,1}}
    \le
    C 2^N  R, 
    \vspace*{3pt}
    \\
    \displaystyle
    \| P_1(t)\|_{M_{3,1}}
    \le 
     \int^t_0 
    \|{\rm div \,} e^{t \mathcal{L}}
    u_{0,N} \|_{M_{3,1}}ds 
    \le 
    C \|u_{0,N}\|_3 
    \le
    C 2^N t R . 
    \vspace*{3pt}
\end{cases}
\]
In the case when $k=2$, it follows from the bilinear estimate that
\[
\begin{aligned}
    \|\Theta_2(t)\|_{M_{3,1}}
    &\le
    C \int^t_0 \|\nabla U_1(s)\|^2_{M_{3,1}} ds 
    \le 
    C 2^{2N} t R^2
    \le C 2^{2N} t R^2,
    \vspace*{3pt}
    \\
    \|U_2(t)\|_{M_{3,1}} 
    &\le 
     \int^t_0 (
    C 2^N\|U_1(s)\|^2_{M_{3,1}}
    + 
    C 2^{N+1} \|\Theta_2(s)\|_{M_{3,1}} 
    + C \|P_1(s)\|_{M_{3,1}} \|\partial_t U_1(s)\|_{M_{3,1}})ds
    \\
    &\le
   C(2^N t R^2 
    +
    2^{3N} t^2 R^2 
    +
     2^{3N} t^2 R^2)
    \le
    C 2^{N} t R^2,
    \\
    \|P_2(t)\|_{M_{3,1}} 
    &\le
    \int^t_0 (
    C 2^{N+1} \|U_2(s)\|_{M_{3,1}} 
    +
    C 2^N\|P_1(s)\|_{M_{3,1}}\|U_1\|_{M_{3,1}} )ds
    \\
    &\le
    C
    (1+C \varepsilon_0)
    2^{2N} t^2 R^2 
    +
    \frac{1}{2}C 2^{2N} t^2R^2
    \le
    C 2^{2N} t^2 R^2. 
\end{aligned}
\]
To prove in the case when $k \ge 3$, we need to estimate $\partial_t U_k, \partial_t \Theta_k$.
By maximal regularity Proposition 2.2, we have the estimate for $t \leq T$ with $T =\varepsilon_0 2^{2N}$, 
\[
\begin{aligned}
    &\|\partial_t \Theta_2 \|_{L^1(0,t;M_{3,1})} +
    \|\Delta \Theta_2 \|_{L^1(0,t;M_{3,1})} 
    \le C 2^{2N} t R^2,
    \\
    &\|\partial_t U_2 \|_{L^1(0,t;M_{3,1})} +
    \|\mathcal{L} U_2 \|_{L^1(0,t;M_{3,1})} 
    \le C 2^{N} t R^2.
\end{aligned}
\]
It is now sufficient to prove the following inequalities to complete the proof of Proposition~\ref{103103}. 

\begin{lem}
\label{103104}
There exist positive constant $c_1, c_2$, such that for sufficiently small $\varepsilon_0 >0$ and $0 < t <T = \varepsilon_0 2^{-2N}$, 
\begin{equation}
    \label{103105}
    \| P_k(t)\|_{M_{3,1}}
    \le
    c_1^{k-1}c_2^k t^k 2^{kN}R^k (1+k)^{-4}, 
\end{equation}
    \begin{equation}
    \begin{aligned}
        \label{103106}
        &\| U_k(t)\|_{M_{3,1}} +
        \|\partial_t U_k(\cdot)\|_{L^1(0,t;M_{3,1})} +
        \| \mathcal{L}U_k(\cdot)\|_{L^1(0,t;M_{3,1})}
        \\
        &\hspace*{11pt}\le
        c_1^{k-1}c_2^{k-1} t^{k-1} 2^{(k-1)N}R^k (1+k)^{-4},   
    \end{aligned}
    \end{equation}
    \begin{equation}
    \begin{aligned}
        \label{103107}
        &\| \Theta_k(t)\|_{M_{3,1}} +
        \|\partial_t \Theta_k(\cdot)\|_{L^1(0,t;M_{3,1})} +
        \| \kappa \Delta \Theta_k(\cdot)\|_{L^1(0,t;M_{3,1})}
        \\
        &\hspace*{11pt}\le
        c_1^{k-1}c_2^k t^{k-1} 2^{kN}R^k (1+k)^{-4}.
    \end{aligned}
    \end{equation}
\end{lem}

\noindent
{\bf Proof of Lemma \ref{103104}. }
We assume the estimate for $1,2,...k-1$ with $k \ge 3$, and show the case of  $k$ to apply an induction argument.
We notice from the definition of $P_k, U_k, \Theta_k$ that $\Theta_k$ can be handled by the assumption for $1,2, ..., k-1$, while $U_k$ needs the estimate of $\Theta_k$, and $P_k$ needs the estimate of $U_k$.

In the proof below, we denote by $c_M$ 
appearing in the bilinear estimate in modulation spaces 
and by $c_2$ the biggest absolute constant 
such that the inequalities before Proposition~\ref{103104} holds. 
We will also apply 
\begin{equation}\label{122205}
 \dfrac{1}{k-1}    \sum_{k_1+k_2 =k} k_1 k_2
     (1+k_1)^{-4}(1+k_2)^{-4} 
\leq \dfrac{C}{(1+k)^4},
\end{equation}
and take $c_1$ such that $c_1 > 10Cc_M$. 

We first show \eqref{103107}. 
\[
\begin{aligned}
    \|\Theta_k(t)\|_{M_{3,1}} &+
        \|\partial_t \Theta_k(t)\|_{L^1(0,t;M_{3,1})} +
        \| \kappa \Delta \Theta_k(t)\|_{L^1(0,t;M_{3,1})}
        \\
        \le&
        c_Mc_2\int^t_0
        \bigg\{
            \sum_{k_1+k_2 =k}
            \|U_{k_1}\|_{M_{3,1}}
            k_2 2^N \|\Theta_{k_2}\|_{M_{3,1}}
            \\
            &+ 
            c_M^2c_2\sum_{k_1+k_2+k_3 =k} 
            \|P_{k_3}\|_{M_{3,1}}
            \|U_{k_1}\|_{M_{3,1}}
            k_2 2^N \|\Theta_{k_2}\|_{M_{3,1}}
            \\
            &+
            c_Mc_2\sum_{k_1+k_2 =k}
            \|\Theta_{k_1}\|_{M_{3,1}}
            k_2 2^N \|U_{k_2}\|_{M_{3,1}}   
            \\
            &+
            c_M^2c_2 \sum_{k_1+k_2+k_3=k} 
            \|P_{k_3}\|_{M_{3,1}}
            \|\Theta_{k_1}\|_{M_{3,1}}
            k_2 2^N \|U_{k_2}\|_{M_{3,1}}    
            \\
            &+
            c_M  \sum_{k_1+k_2 =k}
            \|P_{k_1}\|_{M_{3,1}}
             \|\partial_t\Theta_{k_2}\|_{M_{3,1}}
             \\
            &+ 
            c_M c_2^2  \sum_{k_1+k_2 =k}
            k_1 2^{N} \| U_{k_1}\|_{M_{3,1}}
            k_2 2^{N} \| U_{k_2}\|_{M_{3,1}} 
        \bigg\} ds.
\end{aligned}
\]
We estimate the right hand side above.
The final term is estimated as follows. 
\[
\begin{aligned}
    &\int^t_0 
c_M c_2^2  \sum_{k_1+k_2 =k}
            k_1 2^{N} \| U_{k_1}\|_{M_{3,1}}
            k_2 2^{N} \| U_{k_2}\|_{M_{3,1}}  ds
     \\
     \le&
    \int_0^t s^{k-2}ds \cdot 
     c_M c_2^2 c_1^{k-2} c_2^{k-2} 2^{kN} R^k 
     \sum_{k_1+k_2 =k} k_1 k_2
     (1+k_1)^{-4}(1+k_2)^{-4}
     \\
     \le&
    \dfrac{1}{k-1} c_1^{k-2}
     c_M c_2^{k} 2^{kN} R^k 
     \sum_{k_1+k_2 =k} k_1 k_2
     (1+k_1)^{-4}(1+k_2)^{-4}. 
\end{aligned}
\]
and \eqref{122205} implies that 
\[
\begin{aligned}
    &\int^t_0 
c_M c_2^2  \sum_{k_1+k_2 =k}
            k_1 2^{N} \| U_{k_1}\|_{M_{3,1}}
            k_2 2^{N} \| U_{k_2}\|_{M_{3,1}}  ds
     \\
     \le&
    C c_1^{k-2} c_M c_2^{k} t^{k-1} 2^{kN} R^k 
     (1+k)^{-4}
     = \dfrac{C}{c_1} c_1^{k-1} c_M c_2^{k} t^{k-1} 2^{kN} R^k 
     (1+k)^{-4}.  
\end{aligned}
\]
\color{black}
By estimating others with the bilinear estimate and the assumption of the induction, we have 
\[
\begin{aligned}
    &\|\Theta_k\|_{L^{\infty}(0,t;M_{3,1})} 
    \le
   \dfrac{10Cc_M}{c_1} c_1^{k-1} c_M c_2^{k} t^{k-1} 2^{kN} R^k (1+k)^{-4}. 
\end{aligned}
\]
We know $10Cc_M / c_1 \leq 1$ and conclude the 
estimate for $\Theta_k$. 
We can also show the estimate for $U_k$ in \eqref{103106} in an analogous way to $\Theta_k$. 
The modification of the constant also appears in the estimate of $P_k$. 
We take $c_1>1$ and $c_2>1$ such that $c_1c_2 > 5C(c_1+c_M)$. 
\color{black}
By the assumption of the definition of the induction and applying \eqref{103106}, we see that 
\[
\begin{aligned}
    \|P_k(t)\|_{M_{3,1}} 
    &\le
    \int^t_0 \bigg\|
    {\rm div \,} U_k(s)
    \bigg\|_{M_{3,1}} ds
    + 
    \sum_{k_1+k_2=k}
    \int^t_0 \bigg\|
    {\rm div \,} (P_{k_1}(s)U_{k_2}(s))
    \bigg\|_{M_{3,1}} ds
    \\
    &\le
     Cc_1^{k-1} c_2^{k-1} t^{k} 2^{kN}R^k(1+k)^{-4} 
    \\
    &+ Ck2^N \sum_{k_1+k_2=k} 
    \int^t_0 \bigg\|
     (P_{k_1}(s)U_{k_2}(s))
    \bigg\|_{M_{3,1}} ds.
    \\
    &\le 
    Cc_1^{k-1} c_2^{k-1} t^{k} 2^{kN}R^k(1+k)^{-4} 
    \\
    &+ C k 2^N \int^t_0 s^{k-1}ds  c_1^{k-2}c_M c_2^{k-1} 2^{(k-1)N} R^k \sum_{k_1+k_2=k} (1+k_1)^{-4}(1+k_2)^{-4}
    \\
    &\le
    \frac{C(c_1+c_M)}{c_1c_2}c_1^{k-1} c_2^{k} t^{k} 2^{kN}R^k(1+k)^{-4} 
    \\
    &\le
    c_1^{k-1} c_2^{k} t^{k} 2^{kN}R^k(1+k)^{-4}.
\end{aligned}
\]
We know $5C(c_1+c_M) < c_1c_2$. 
Then we showed the inequality \eqref{103105}.

\vskip3mm 

\noindent
{\bf Proof of Theorem 1.1} 
For $t \le \varepsilon_0 2^{-2N}$, 
We have
\begin{equation}
\begin{aligned}
    \label{110101}
    \|\theta(t)\|_{\dot B^{-1}_{3,1}}
    &\ge
    \bigg\|
    \sum_{ -\delta N \le j \le 0} \phi_j * \theta (t)
    \bigg\|_{\dot B^{-1}_{3,1}}
    \\
    &\ge
    \bigg\|
    \sum_{ -\delta N \le j \le 0} \phi_j * \Theta_2 (t)
    \bigg\|_{\dot B^{-1}_{3,1}} - 
    \sum_{k \ge 3}
    \bigg\|
    \sum_{j=-\delta N}^0 \phi_j * \Theta_k (t)
    \bigg\|_{\dot B^{-1}_{3,1}}
\end{aligned}
\end{equation}
The first term is handled by \eqref{081401}, and it follows from the embedding theorem that
\[
\begin{aligned}
    \bigg\|
    \sum_{ -\delta N \le j \le 0} \phi_j * \Theta_k (t)
    \bigg\|_{\dot B^{-1}_{3,1}}
    &\le
    \sum_{ -\delta N \le j \le 0}
    \bigg\|
     \phi_j * \Theta_k (t)
    \bigg\|_{\dot B^{-1}_{3,1}}
    \le
    \sum_{ -\delta N \le j \le 0} \sum_{l \in \mathbb{Z}} 
    2^{-l}\|\phi_l *\phi_j * \Theta_k(t)\|_3
    \\
    &\le
    C \sum_{ -\delta N \le j \le 0}
    2^{-j}\|\phi_j * \Theta_k(t)\|_3
    \le
    C \sum_{ -\delta N \le j \le 0}
    2^{-j}\|\Theta_k(t)\|_{M_{3,1}}
    \\
    &\le C c_1^{k-1}c_2^k t^{k-1} 2^{kN}R^k (1+k)^{-4}
    \sum_{ -\delta N \le j \le 0} 2^{-j}
    \\
    &\le
    C c_1^{k-1}c_2^k \varepsilon_0^{k-1} 2^{(-k+2+\delta)N} R^k (1+k)^{-4}. 
\end{aligned}
\]
Using above inequality to \eqref{110101}, we have
\[
\begin{aligned}
    \|\theta(t)\|_{\dot B^{-1}_{3,1}}
    &\ge
    c\delta R^2 N^{\frac{1}{3}} -
    (R N^{-\frac{1}{3}})^2 2^{-cN} -
    C\sum_{k\ge 3} c_1^{k-1}c_2^k \varepsilon_0^{k-1} 2^{(-k+2+\delta)N} R^k (1+k)^{-4}
    \\
    &\ge
    c\delta R^2 N^{\frac{1}{3}} -
    (R N^{-\frac{1}{3}})^2 2^{-cN} -
    C(\varepsilon_0 c_1)^2 c_2^3 R^3 2^{(\delta -1)N} 
    \sum_{k\ge 0} (2^{-N}c_1c_2\varepsilon_0
    R)^k
    \\
    &\to \infty \qquad as \ 
    N \to \infty.
\end{aligned}
\]
\qed

\appendix 

\section{The bound from below in general dimensions}

In this section, we prove the corresponding estimate  to \eqref{110402} for the $d$-dimensional case to clarify the dependence on the dyadic frequency parameter.

\begin{lem}
  \label{0125-01}
  \begin{equation}
  \begin{aligned}
    \label{122102}
    \sum_{ -\delta N \le j \le 0} 
    \sum_{k \neq j }  
    2^{-j} 2^{-(d-1)j} 2^{-(d-1)k}
      & \bigg\|
        \phi_j * 
        \bigg(
         \phi_j (x - 2^{|j| +2N}e_1)  
         \phi_k (x - 2^{|k| +2N}e_1)  {\rm sin \,}^2(2^Nx_1) 
         \bigg)
         \bigg\|_d\\
  &\leq C N 2^{-d(1- \delta)N}. 
  \end{aligned}
  \end{equation}
  \begin{equation}
    \begin{aligned}
      \label{122104}
      \sum_{ -\delta N \le j \le 0} 
      \sum_{k \neq j} 
      \sum_{l \neq j} 
      2^{-j} 2^{-(d-1)l} 2^{-(d-1)k}
         & \bigg\|
          \phi_j * 
          \bigg(
           \phi_l (x - 2^{|l| +2N}e_1)  
           \phi_k (x - 2^{|k| +2N}e_1)  {\rm sin \,}^2(2^Nx_1) 
           \bigg)
           \bigg\|_d\\
    &\leq C N^2 2^{-d(1- \delta)N}.
    \end{aligned}
    \end{equation}
  \begin{equation}
  \begin{aligned}
    \label{122103}
    \sum_{ -\delta N \le j \le 0}\sum_{k \neq j} 
    2^{-j} 2^{-2(d-1)k}
    &\bigg\|\phi_j *   
    \bigg|(\phi_k(x- 2^{|k|+2N}e_1))\bigg|^2 {\rm sin \,}^2(2^Nx_1)\bigg\|_{L^d(A_j)}
    \\
  &\leq C2^{-d(1-2\delta)N},
  \end{aligned}
  \end{equation}
  where we define the set $A_j$ by
  \begin{equation}
    \label{0124-01}
    A_j := 
    \{
      |x - 2^{|j|+2N}e_1| \le 2^{-j}
    \}
  \end{equation}
  \end{lem}
\color{black}
  
  \begin{proof}
  We first show \eqref{122102}.
  By Young's inequality, we have
  \[
  \begin{aligned}
    \bigg\|
    \phi_j * 
    \bigg(
     \phi_j (x - 2^{|j| +2N}e_1)  
     \phi_k (x - 2^{|k| +2N}e_1)  {\rm sin \,}^2(2^Nx_1) 
     \bigg)
     \bigg\|_d
     \\
     \le
     c\|
     \phi_j (x - 2^{|j| +2N}e_1)  
     \phi_k (x - 2^{|k| +2N}e_1)
     \|_d. 
  \end{aligned}
  \]
  If $x\in B_{2^N}(2^{|j|+2N}e_1)$, it satisfies 
  \[
    2^k|x - 2^{|k| +2N} e_1|\ge 
    c2^k|2^{|j|+2N}- 2^{|k|+2N}e_1| =2^k2^{2N},
  \] 
  and there exists some $\alpha >d$ such that
  \begin{equation}
    \label{122702}
    |\phi_k (x - 2^{|k| +2N}e_1)|
    \le 2^{kd}(1+2^k 2^{2N})^{-\alpha}
    \le 2^{-2dN}. 
  \end{equation}
  If $x\in B^c_{2^N}(2^{|j|+2N}e_1)$, the analogous argument leads to
  \begin{equation}
    \label{122701}
    |\phi_j (x - 2^{|j| +2N}e_1)|
    \le 2^{jd}(1+2^j 2^{N})^{-\alpha}
    \le 2^{-dN}. 
  \end{equation}
  Since $-\delta N \le j,k \le 0$,  we have
  \[
  \begin{aligned}
     &\|
     \phi_j (x - 2^{|j| +2N}e_1)  
     \phi_k (x - 2^{|k| +2N}e_1)
     \|_d
     \\
     \le
     &\|\phi_j (x - 2^{|j| +2N}e_1)  
     \phi_k (x - 2^{|k| +2N}e_1)
     \|_{L^d(B_{2^N}(2^{|j|+2N}e_1))}
     \\
     &+\|\phi_j (x - 2^{|j| +2N}e_1)  
     \phi_k (x - 2^{|k| +2N}e_1)
     \|_{L^d(B^c_{2^N}(2^{|j|+2N}e_1))}
     \\
     \le
     &2^{-2dN} \| \phi_j (x - 2^{|j| +2N}e_1)\|_d
     +
     2^{-dN} \| \phi_k (x - 2^{|k| +2N}e_1)\|_d
     \\
     \le
     &2^{-2dN} 2^{(d-1)j} +2^{-dN}2^{(d-1)k}
     \le C 2^{-dN}2^{(d-1)k}
  \end{aligned}
  \]
  We apply above to the right hand side of \eqref{122102} and then 
  \[
  \begin{aligned}
    \sum_{ -\delta N \le j \le 0} 
    \sum_{k \neq j }  
    &2^{-j} 2^{-(d-1)j} 2^{-(d-1)k}
       \bigg\|
        \phi_j * 
        \bigg(
         \phi_j (x - 2^{|j| +2N}e_1)  
         \phi_k (x - 2^{|k| +2N}e_1)  {\rm sin \,}^2(2^Nx_1) 
         \bigg)
         \bigg\|_d\\ 
         \le
         C\sum_{ -\delta N \le j \le 0} 
         \sum_{k \neq j }  
         &2^{-j} 2^{-(d-1)j} 2^{-dN}
        \le CN 2^{-d(1-\delta)N}, 
  \end{aligned}
  \]
  which proves \eqref{122102}.
  The estimate \eqref{122104} also follows from a similar argument. 
  Finally we show the estimate \eqref{122103}.
  We separate \eqref{122103} To
  \begin{equation}
    \begin{aligned}
      \label{012401}
      &\bigg\|\phi_j *   
      \bigg|(\phi_k(x- 2^{|k|+2N}e_1))\bigg|^2 \bigg\|_{L^d(A_j)} 
      \\
      \le &
      \bigg\| \int_{B_{2^N}(2^{|j|+2N}e_1)}
      \phi_j (y)(\phi_k(x-y- 2^{|k|+2N}e_1))^2 dy 
      \bigg\|_{L^d(A_j)}
      \\
      &+
      \bigg\| \int_{B^c_{2^N}(2^{|j|+2N}e_1)}
      \phi_j (y)(\phi_k(x-y- 2^{|k|+2N}e_1))^2 dy 
      \bigg\|_{L^d(A_j)} 
    \end{aligned}
  \end{equation}
  where $A_j$ is defined by \eqref{0124-01}.
  \color{black}
  If $y\in B_{2^N}(2^{|j|+2N}e_1)$, there exists some $\alpha >d$ such that $\phi_k$ satisfies \eqref{122702} and the first term of the right hand side above satisfies
  \[
  \begin{aligned}
    &\bigg\| \int_{B_{2^N}(2^{|j|+2N}e_1)}
    \phi_j (y)(\phi_k(x-y- 2^{|k|+2N}e_1))^2 dy 
    \bigg\|_{L^d(A_j)}
    \\
    \le
    &
    C2^{jd}2^{2kd}
    \bigg\| \int_{B_{2^N}(2^{|j|+2N}e_1)}
     (1+ 2^j|y|)^{-\alpha} 
     (1+ 2^k|x-y|)^{-2\alpha}dy 
     \bigg\|_{L^d(A_j)}
    \\
    \le
    &
    C2^{jd}2^{2kd}
    \bigg\| \int_{B_{2^N}(2^{|j|+2N}e_1)}
     (1+ 2^j|y|)^{-\alpha} 
     (1+ 2^k2^{2N})^{-2\alpha}dy 
     \bigg\|_{L^d(A_j)}
    \\
    \le
    &C 2^{jd}2^{2kd}
    \bigg\|\frac{2^{-jd}}{(2^k 2^{2N})^{2d}}\bigg\|_{L^d(A_j)}
    \le 
    C 2^{jd}2^{2kd}2^{-jd}2^{-kd}2^{-2dN} 2^{-j}
    =
    C2^{-2dN} 2^{-j}2^{kd}.
  \end{aligned}
  \]
  We obtain the upper bound for the first term of the right hand side on \eqref{012401}
  \[
  \begin{aligned}
    &\sum_{ -\delta N \le j \le 0}\sum_{k \neq j} 
    2^{-j} 2^{-2(d-1)k}
    \bigg\| \int_{B_{2^N}(2^{|j|+2N}e_1)}
    \phi_j (y)(\phi_k(x-y- 2^{|k|+2N}e_1))^2 dy 
    \bigg\|_{L^d(A_j)}
    \\
    \le
    &C 2^{-2dN}\sum_{ -\delta N \le j \le 0}\sum_{k \neq j} 
    2^{-j} 2^{-2(d-1)k}
    2^{-j}2^{kd}
    \le
    C 2^{-2dN}\sum_{ -\delta N \le j \le 0}\sum_{k \neq j} 
    2^{-2j}2^{-(d-2)k}
    \\
    \le&
    C2^{-d(2-\delta)N}.
  \end{aligned}
  \]
  If $y\in B^c_{2^N}(2^{|j|+2N}e_1)$, there also exists some $\alpha >d$ such that $\phi_j$ satisfies \eqref{122701} and the second term of the right hand side of \eqref{012401} satisfies
  \[
  \begin{aligned}
    &\bigg\| \int_{B^c_{2^N}(2^{|j|+2N}e_1)}
    \phi_j (y)(\phi_k(x-y- 2^{|k|+2N}e_1))^2 dy 
    \bigg\|_{L^d(A_j)} 
    \le 
    C 2^{jd}2^{2kd}
    \bigg\|\frac{2^{-2kd}}{(2^j2^{N})^d}\bigg\|_{L^d(A_j)}
    \\
    \le
    &C2^{-dN} 2^{-j}.
  \end{aligned}
  \]
  We also obtain the upper bound for second term of the right hand side on 
  \eqref{012401}
  \[
  \begin{aligned}
    &\sum_{ -\delta N \le j \le 0}\sum_{k \neq j} 
    2^{-j} 2^{-2(d-1)k}
    \bigg\| \int_{B^c_{2^N}(2^{|j|+2N}e_1)}
    \phi_j (y)(\phi_k(x-y- 2^{|k|+2N}e_1))^2 dy 
    \bigg\|_{L^d(A_j)} 
    \\
    \le 
    &C 2^{-dN} \sum_{ -\delta N \le j \le 0}\sum_{k \neq j}  2^{-2j} 2^{-2(d-1)k}
    \le
    C2^{-d(1-2\delta)N}.
  \end{aligned}
  \]
The proof of Lemma~\ref{0125-01} is completed.
  \end{proof}

\begin{prop}
 Let $d \ge 2$. 
 For the initial data $u_{0,N} \in (\mathcal S(\mathbb R^d))^d$ such that
  \[
    u_{0,N} (x) :=
    \bigg(R N^{-\frac{1}{d}} \sum_{ -\delta N \le j \le 0}
    2^{-(d -1)j} \phi_j (x - 2^{|j| +2N}e_1) {\rm sin} (2^Nx_1) ,0, \cdots, 0 \bigg), 
\]
there exist $c, \varepsilon_0 > 0$ such that 
\begin{equation}
  \label{081401} 
  \bigg\| \int^t_0 e^{(t - s)\mathcal L} | \nabla 
  e^{ s\mathcal L} u_{0,N} |^2 ds 
   \bigg\|_{\dot B^{-1}_{d,1}(\mathbb{R}^d)} \ge 
   c\delta R^2 N^{\frac{1}{d}} -
   (R N^{-\frac{1}{d}})^2 2^{-(1-2\delta)N}, 
  \end{equation}
for $t = \varepsilon_0 2^{-2N}$. 
\end{prop}

\noindent
 \begin{proof}
We start by estimate the term from the Gaussian. 
We write 
  \[
  \begin{aligned}
  &\int^t_0 e^{(t - s)\mathcal L} | \nabla 
  e^{ s\mathcal L} u_{0,N} |^2 ds 
  \\
  =
  &
  \mathcal{F}^{-1}
  \Big[
    \int_{\mathbb R^d}
    \widehat{\varphi} (\xi) 
    i(\xi-\eta)i \eta
  \widehat u_{0,N}(\xi-\eta)
  \widehat u_{0,N}(\eta) e^{-t|\xi|^2}
  e^{-t|\xi|^2}
  \frac{1-e^{-t(|\xi - \eta|^2+|\eta|^2 -|\xi|^2)}}{|\xi - \eta|^2+|\eta|^2 -|\xi|^2}d\eta\Big](x).
  \end{aligned}
  \]
  For $\varepsilon_0 >0$ to be taken later, 
  we apply the Talor expansion with 
  taking $t = \varepsilon_0 2^{-2N}, |\xi-\eta|, |\eta| \simeq 2^N$ and $|\xi| \leq 1$, 
and we will use 
  \color{black}
\begin{equation}\label{122201}
    \begin{split}
      &\displaystyle
      i(\xi - \eta) \cdot i\eta
      \frac{1-e^{-t(|\xi - \eta|^2+|\eta|^2 -|\xi|^2)}}{|\xi - \eta|^2+|\eta|^2 -|\xi|^2}
     \\ 
     &= \displaystyle  
      \dfrac{|\xi - \eta|^2+|\eta|^2 -|\xi|^2}{2}       \sum_{m \in \mathbb{N}_{\ge 0}}
      \frac{(-1)^{m-1}}{m!}t^m
      (|\xi - \eta|^2+|\eta|^2 -|\xi|^2)^{m-1}
      \\
      &
      =
      \displaystyle
      \frac{1}{2}
      \sum_{m \in \mathbb{N}}
      \frac{(-1)^{m-1}}{m!}t^m
      (|\xi - \eta|^2+|\eta|^2 -|\xi|^2)^m
      \\
      &
      \displaystyle
      = \frac{1}{2}
      \bigg(\varepsilon_0 - 
      \sum_{m\ge2} \frac{(-1)^{m}}{m!}
      O(C^m \varepsilon_0^m )
      \bigg).
\end{split}
\end{equation}

We show \eqref{081401}. 
  Let $\varphi \in \mathcal{S}(\mathbb{R}^3)$ satisfy $\supp \widehat{\varphi} \subset \{\xi \in \mathbb{R}^d||\xi| \le 1\}$, 
  consider the restriction to low frequency part and estimate the left hand side of \eqref{081401} from below. 
  \begin{equation}
    \begin{aligned}
      \label{110801}
      & \bigg\| \int^t_0 e^{(t - s)\mathcal L} | \nabla 
      e^{ s\mathcal L} u_{0,N} |^2 ds 
       \bigg\|_{\dot B^{-1}_{d,1}(\mathbb{R}^d)}
       \ge 
       c \bigg\|\varphi * \int^t_0 e^{(t - s)\mathcal L} | \nabla 
      e^{ s\mathcal L} u_{0,N} |^2 ds 
       \bigg\|_{\dot B^{-1}_{d,1}(\mathbb{R}^d)}
       \\
       \ge 
& c      \sum_{ -\delta N \le j \le 0}2^{-j}
       \bigg\|
       \phi_j* \varphi * \int^t_0 e^{(t - s)\mathcal L} | \nabla 
      e^{ s\mathcal L} u_{0,N} |^2 ds
       \bigg\|_d
       \\
       =
       & (R N^{-\frac{1}{d}})^2 
       \sum_{-\delta N \leq j \leq 0}^0 2^{-j}
       \bigg\|
       \phi_j * \varphi *\int^t_0
       e^{(t-s)\mathcal{L}}
       |\nabla e^{s\mathcal{L}}
       \sum^0_{k = -\delta N } 
       2^{-(d-1)k} \phi_k (x - 2^{|k| +2N}e_1)  {\rm sin \,}(2^Nx_1) |^2ds
        \bigg\|_d
\end{aligned}
\end{equation}
Here, we divide the sum over $- \delta N \leq k \leq 0$ into three cases with an elementary equality that 
for each $j$
\[
\Big| \sum_{-\delta N \leq k \leq 0} a_k \Big|^2 
= a_j ^2 + 2 a_j \sum_{k \not = j} a_k 
+ \Big| \sum_{k \not = j} a_k \Big|^2.  
\]
  We first consider the third case above. 
We apply the Fourier multiplier theorem to 
\[
\varphi(\xi)  \cdot       i(\xi - \eta) \cdot i\eta e^{-t |\xi|^2}
   \frac{1-e^{-t(|\xi - \eta|^2+|\eta|^2 -|\xi|^2)}}{|\xi - \eta|^2+|\eta|^2 -|\xi|^2}, 
\]
and then we have the bound from above with the 
leading term $\varepsilon _0$, roughly due to the expansion \eqref{122201}. Therefore, we obtain 

  \begin{equation}
    \begin{aligned}
      \label{110802}
      & (R N^{-\frac{1}{d}})^2 
     \sum_{ -\delta N \le j \le 0}2^{-j}
      \bigg\| \phi_j * 
      \mathcal{F}^{-1}
      \bigg[
        \widehat{\varphi}(\xi)
       (\eta -\xi) \eta e^{-t|\xi|^2}
       \frac{(1-e^{-t(|\xi - \eta|^2+|\eta|^2 -|\xi|^2)})}{|\xi - \eta|^2+|\eta|^2 -|\xi|^2}
      \\
      &\hspace{120pt}
       \mathcal{F}
       \bigg[
       \bigg|\sum^0_{k \neq j} 
       2^{-(d-1)k} \phi_k (x - 2^{|k| +2N}e_1)  {\rm sin \,}(2^Nx_1) 
       \bigg|^2
      \bigg](\xi) 
      \bigg]
      \bigg\|_d
       \\
       \le&
       C \varepsilon_0(R N^{-\frac{1}{d}})^2 
     \sum_{ -\delta N \le j \le 0}2^{-j}
      \bigg\|
      \phi_j * 
      \bigg|
      \sum^0_{k \neq j } 2^{-(d-1)k}
       \phi_k (x - 2^{|k| +2N}e_1)  {\rm sin \,}(2^Nx_1) 
       \bigg|^2
       \bigg\|_d
\\
 \leq &      c\varepsilon_0R^2 N^{\frac{1}{d}} 2^{-d(1- 2\delta)N},
    \end{aligned}
  \end{equation}
where we have applied \eqref{122103} in the last inequality. 
  We also obtain the upper bound for the second case  
  with the estimate \eqref{122102} or 
  \eqref{122104}. 

\color{black}
Finally we estimate the first case from below, 
and need to consider the following. 
\[
F=          \frac{1}{2}\varepsilon_0
         (R N^{-\frac{1}{d}})^2 
        \sum_{ -\delta N \le j \le 0}2^{-j} 2^{-2(d-1)j} 
         \bigg\|
         \phi_j * \varphi * 
          \phi_j^2(x-2^{|j|+2N}e_1)
          {\rm sin \,}^2 (2^Nx_1)
         \bigg\|_d.
\]
Let us recall 
\[
  A_j := 
  \{
    |x - 2^{|j|+2N}e_1| \le 2^{-j}
  \},
\]
\color{black}
  and show the estimate $F$.
  We restrict the $L^d$ norm in $A_j$ which written in \eqref{0124-01}, 
  \begin{equation}
  \begin{aligned}
    \label{102301}
      F =&
      \frac{\varepsilon_0}{2}(
        RN^{-\frac{1}{d}})^2
     \sum_{ -\delta N \le j \le 0}2^{-j}   2^{-2(d-1)j}
      \|\phi_j * \varphi*(\phi_j^2(x- 2^{|j|+2N}e_1)) {\rm sin \,}^2(2^Nx_1)\|_d 
      \\
      \ge &
      \frac{\varepsilon_0}{2}
      (RN^{-\frac{1}{d}})^2
     \sum_{ -\delta N \le j \le 0}2^{-j}   2^{-2(d-1)j}
      \|\phi_j *(\phi_j^2(x- 2^{|j|+2N}e_1)) {\rm sin \,}^2(2^Nx_1)\|_{L^d(A_j)}
      \\
      \ge
      &
      \frac{\varepsilon_0}{2}
      (RN^{-\frac{1}{d}})^2
     \sum_{ -\delta N \le j \le 0}2^{-j}   2^{(d+2)j}
       \bigg\{
      \bigg\|\int_{\mathbb R^d} \phi_0 (2^j(x-y))
      \Big(\phi_0(2^j(y- 2^{|j| +2N}))\Big)^2 dy \bigg\|_{L^d(A_j)}\\
      & \hspace{30pt}
      - 
      \frac{1}{2}\bigg\| \int_{\mathbb R^d} \phi_0 (2^j(x-y))
      \Big(\phi_0(2^j(y- 2^{|j| +2N}))\Big)^2 {\rm cos \,}(2^{N+1}y_1)
      dy
      \bigg\|_{L^d(A_j)}
      \bigg\}.
  \end{aligned}
  \end{equation}
  We write the first term on the right hand side of \eqref{102301}
  \[
  \begin{aligned}
      &\frac{\varepsilon_0}{2}
      (RN^{-\frac{1}{d}})^2
     \sum_{ -\delta N \le j \le 0}2^{-j}   2^{(d+2)j}
      \bigg\|\int_{\mathbb R^d} \phi_0 (2^j(x-y))
      \Big(\phi_0(2^j(y- 2^{|j| +2N}))\Big)^2 dy \bigg\|_{L^d(A_j)}
      \\
      =&\frac{\varepsilon_0}{2}
      (RN^{-\frac{1}{d}})^2
     \sum_{ -\delta N \le j \le 0}2^{-j}   2^{(d+2)j}  
        2^{-(d+1)j}
        \bigg\|\int_{\mathbb R^d} \phi_0 (x-y)\Big(\phi_0(y)\Big)^2 dy \bigg\|_{L^d(\{|x|\le1\})} 
      \\
      =&
      \frac{\varepsilon_0}{2}
      (RN^{-\frac{1}{d}})^2
     \sum_{ -\delta N \le j \le 0}
        \bigg\|\int_{\mathbb R^d} \phi_0 (x-y)\Big(\phi_0(y)\Big)^2 dy \bigg\|_{L^d(\{|x|\le1\})}. 
  \end{aligned}
  \]
  The second term on the right hand side of \eqref{102301} is dealt with the integration by parts, we have that
  \[
  \begin{aligned}
    & \frac{1}{2}
    (RN^{-\frac{1}{d}})^2
   \sum_{ -\delta N \le j \le 0}2^{-j}   2^{(d+2)j}\bigg\| \int_{\mathbb R^d} \phi_0 (2^j(x-y))\Big(\phi_0(2^j(y- 2^{|j| +2N}))\Big)^2 {\rm cos \,}(2^{N+1}y_1)
    dy
    \bigg\|_{L^d(A_j)}\\
    =& \frac{1}{2}
    (RN^{-\frac{1}{d}})^2
   \sum_{ -\delta N \le j \le 0}2^{-j}   2^{(d+2)j}
      \bigg\| \int_{\mathbb R^d} \partial_{y_1}
      \bigg(\phi_0 (2^j(x-y))\Big(\phi_0(2^j(y- 2^{|j| +2N}))\Big)^2 \bigg)
      \frac{{\rm sin \,}(2^{N+1}y_1)}{2^{N+1}}
      dy
      \bigg\|_{L^d(A_j)}
      \\
      =& \frac{1}{4}
      (RN^{-\frac{1}{d}})^2
     \sum_{ -\delta N \le j \le 0}2^{-j}   2^{(d+2)j}
       2^j2^{-N}    
       \bigg\| \int_{\mathbb R^d} 
      \bigg(\phi_0 (2^j(x-y))\Big(\phi_0(2^j(y- 2^{|j| +2N}))\Big)^2 \bigg) {\rm sin \,}(2^{N+1}y_1)
      dy
      \bigg\|_{L^d(A_j)}\\
      \le&\frac{1}{4}
      (RN^{-\frac{1}{d}})^2
     \sum_{ -\delta N \le j \le 0}2^{-j}   2^{(d+2)j}
       2^j2^{-N} 2^{-(d+1)j}   
      \bigg\| \int_{\mathbb R^d} 
      \bigg(\phi_0 (x-y)\Big(\phi_0(y)\Big)^2 \bigg) 
      dy
      \bigg\|_{L^d(\{|x|\le1\})}
      \\
      =&
      \frac{1}{4}
      (RN^{-\frac{1}{d}})^2
     \sum_{ -\delta N \le j \le 0}2^{j-N}  
      \bigg\| \int_{\mathbb R^d} 
      \bigg(\phi_0 (x-y)\Big(\phi_0(y)\Big)^2 \bigg) 
      dy
      \bigg\|_{L^d(\{|x|\le1\})}
  \end{aligned}
  \]

  Then these inequalities above imply that
  \begin{equation}
  \begin{aligned}
  \label{102401}
     F
     &\ge
  \frac{c}{4}R^2 N^{-\frac{2}{d}}
  \varepsilon_0 
  \sum_{ -\delta N \le j \le 0} (1- 2^{-1}2^j 2^{-N})
  =\frac{c}{4}R^2 \delta N^{1-\frac{2}{d}}.
  \end{aligned}
  \end{equation}

  Therefore we obtain \eqref{081401} by
  \eqref{110802} and 
  \eqref{102401}.
\end{proof}

\end{document}